\documentclass{amsart}
  \usepackage{paralist}
\usepackage{graphicx}

\newcommand{\N}{\mathbb{N}}
\newcommand{\R}{\mathbb{R}}

\newcommand{\be}{\begin{equation}}
\newcommand{\ee}{\end{equation}}

\newtheorem{theorem}{Theorem}[section]
\newtheorem{proposition}[theorem]{Proposition}

\newtheorem{definition}[theorem]{Definition}

\begin{document}

\thispagestyle{empty}

\title[Forward-backward parabolic problems]{Non-uniqueness results for entropy two-phase solutions of forward-backward parabolic problems with unstable phase}

\author{Andrea Terracina}
\maketitle

\begin{abstract} 
This paper study the  well--posedness of the entropy formulation given by Plotnikov in [\textit{Differential Equations}, 30 (1994), pp.\ 614--622]  for  forward--backward parabolic problem obtained as singular limit of a proper  pseudoparabolic approximation. It was proved in [C.\ Mascia, A.\ Terracina, and A.\ Tesei, \textit{Arch.\ Ration.\ Mech.\ Anal.}, 194 (2009), pp.\ 887--925] that such formulation
gives uniqueness when the solution takes values in the stable phases. Here we consider the situation in which unstable phase is taken in account, proving that, in general,    uniqueness does not hold. 
\end{abstract}

\keywords{phase transition, forward--backward equations, ill--posed problems}

\section{Introduction}\label{intro}

 In this paper we consider  the following    {\it forward--backward}
parabolic problem: 
\be\label{B1}
 \left\{\begin{array}{lll} &u_t= \phi(u)_{xx} & \hbox{ in } Q_T:=\Omega\times (0,T)\\
& \phi(u)_x(0,t)\equiv \phi(u)_x(L,t)\equiv 0 & \hbox{ in } (0,T)\\
&  u(x,0)=u_0(x) & \hbox{ in } \Omega:=(0,L)\end{array}
  \right.
 \ee
where   $\phi$ is a {\it nonmonotone}   function. Obviously this kind of problem is ill--posed
whenever $u$ takes values in the interval in which $\phi$ decreases.   
In particular in this paper we consider a piecewise linear function $\phi$, 
 namely:
\begin{equation}\label{philin}
\phi(u) = \left\{\begin{aligned}
&\phi_1(u) \qquad & \textrm{for}\quad  & u \le b \\
&\phi_0(u)\qquad   & \textrm{for}\quad  & b < u < c\\
&\phi_2(u)\qquad  & \textrm{for}\quad   & u \ge d \, ,
\end{aligned}\right.
\end{equation}
where
\begin{equation*}
\phi_i(u):=\alpha_i\,u+\gamma_i \, ,\,i=1,2, \qquad
\phi_0(u):= \frac{A(u-b) -B(u-c)}{c-b}\; .
\end{equation*}
\begin{figure}[htbp]\label{funzione2}
\begin{center}
\includegraphics[width=9cm, height=5cm]{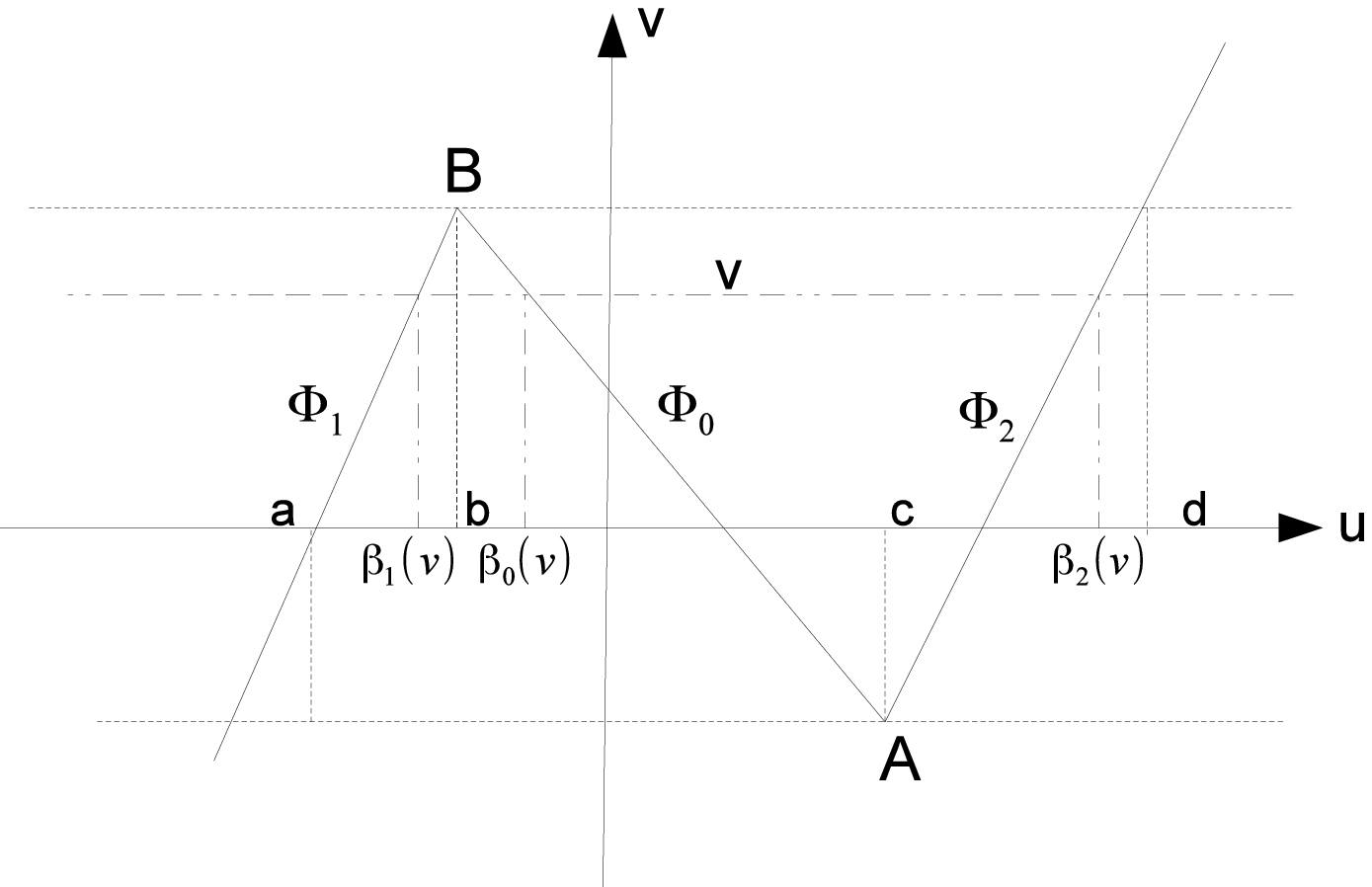}
\caption{The function  $\phi$. }
\end{center}
\end{figure}

Here $-\infty<b<c<\infty$,  $\alpha_i>0$, 
$\gamma_i\in\R$, $i=1,2$,
$A:=\phi_2(c)<\phi_1(b)=:B$.

Let us denote with $\beta_1:(-\infty, B)\longrightarrow \R$, $\beta_2:(A, \infty)\longrightarrow \R$,  $\beta_0:(A, B)\longrightarrow \R$, respectively, the inverse function of $\phi_1$, $\phi_2$ and $\phi_0$ (see Fig.\ref{funzione2}).

\smallskip

The differential equation in  (\ref{B1}) with the response function $\phi$ of cubic type  arises in the theory of phase transition. The function $u$ gives the phase fields and its values characterize the different phases;  
the  half-lines $ (-\infty,b)$ and  $(c,\infty)$ correspond to stable  phases 
and the  interval $(b,c)$ to the unstable one ($e.g.$, see \cite{BS}). 

\noindent In \cite {H} (see also \cite{HN}) it was proved that the problem (\ref{B1}) with a  piecewise function $\phi$ has an infinite number of solutions.  It is interesting to underline that  the solutions given in \cite{H} takes values only in the two stable phases when $t>0$. In some sense the presence of the unstable phase allow  to pass from one phase to the other with too much freedom.

It is important to recall that for  forward--backward  parabolic problems it is possible to state   uniqueness results assuming  that the solution is quite regular (see \cite  {Fr} for the backward case and \cite {Lai1}, \cite {Lai2} for the forward--backward case).

In order to give a good formulation for the problem (\ref{B1}) a natural approach is to introduce a proper regularization, obviously the choice of the regularization terms is related to the physical phenomenon  that we want to describe. In the case of the model of phase transition the original problem is very complicate from a mathematical point of view since there are many terms to take in account. As a matter of fact, it is possible to choose different type of regularizations  in which only some  phenomena are considered (see e.g.  \cite {BFJ}, \cite {BS}, \cite{fife}, \cite{gurtin} \cite{BBDU}, \cite{BFG},   \cite{GG}, \cite{V}). Using this point of view the choice of a particular regularization depends on the phenomena which we want to highlight.

Here we  consider the following pseudoparabolic regularization 
\begin{equation} \label{E2}
 \left\{\begin{array}{lll} &u_t= v_{xx}  & \hbox{ in } Q_T:=(0,L)\times (0,T)\\
& v _x(0,t)\equiv v_x(L,t)\equiv 0 & \hbox{ in } (0,T)\\
&  u(x,0)=u_0(x) & \hbox{ in } (0,L),\end{array}
  \right.\ee
where $v=\phi(u)+\epsilon u_t$, $\epsilon>0$.

\noindent The third order term  in the right hand   side of the differential equation in (\ref{E2})  is a viscosity term  related to   nonequilibium effects (see e.g. \cite{BFJ}, \cite{FJ}, \cite{gurtin}). Let us observe, that in  the regularization (\ref{E2}) it is not considered the characteristic term  of the Cahn--Hilliard equation that describes the cost of the  inhomogeneities in phase transition models. 

It is worth to note that the approximation equations (\ref{E2}) have an independent interest that  is beyond the physical model. More precisely,  when $\phi$ is a linear function these equations suggest a   Yoshida approximation of the differential equation in (\ref{B1}). In fact, this kind of equations were introduced in quasi--reversibility  methods to approximate backward parabolic problem, see \cite {LL},  \cite {ST1}, \cite {E}.

Problem (\ref{E2}) with a general nonlinear function $\phi$ of cubic type  was studied in \cite{NP}, whereas  the singular limit    was  analyzed  by Plotnikov (see \cite{Pl1}, \cite{Pl2}  \cite{Pl3} ), a similar analysis for other type of nonlinearity was considered in  \cite{Sma}. The idea, in analogy with conservation laws, is  to give an entropy formulation of the ill-posed problem (\ref{B1}) assuming that the physical solutions  of it  are that obtained when $\epsilon$ goes to $0^+$ as limit of solutions of problems (\ref{E2}).

In general,  we point out  that functions $u$ obtained as limit of the problems (\ref{E2}) do not   satisfy  equation in  (\ref{B1}) in the classical sense, more precisely we have a solution $(u,\lambda_i,v)$, $i=0,1,2$,  such that
\be \label {lambda}u=\lambda_1\beta_1(v)+\lambda_0\beta_0(v)+\lambda_2\beta_2(v),\ee
where $\lambda_i(x,t)\ge 0$,  $i=0,1,2$, $\displaystyle{\sum_{i=0}^2\lambda_i(x,t)=1}$ in $Q_T$ and the equation $u_t=v_{xx}$ is satisfied in the weak sense.

\noindent In this context the solution $u$ can  be regarded as a superposition of different states and fulfills the differential equation (\ref{B1}) in the sense of the Young measures (see e.g. \cite {Pl3}, \cite {MTT}). 

\noindent  We shall give the precise definition of entropy solution
suggested by Plotnikov in Section 2.   However, we anticipate that this definition provides a good formulation for solutions that takes values  in the two stable phases. In particular, in \cite {arma} \cite {T}, \cite {ST}, it was introduced the  ``two--phase problem''    for which   the initial data   and the solution takes values  in   the two stable phases
and there is a regular  interface that separates different phases. In this situation the entropy formulation  of Plotnikov suggest an admissibility condition for the evolution of the interface. For this type of solutions local existence and uniqueness  with  the response function (\ref{philin}) was obtained in \cite {arma}, global existence was proved in \cite {T}, while   existence, uniqueness and the study of the singular limit  for a general nonlinear cubic type function $\phi$ is established in \cite{ST}.
  It is easy  to check  that  the admissibility condition along the interface does not allow to consider the solutions given in \cite {H}, we have a stricter condition for jumping from one phase to the other and  this guarantees uniqueness.
    
\noindent Using these considerations,  we can guess that the formulation of entropy solution given by Plotnikov could be satisfactory also in the general case where   unstable phase is taken in account. It is necessary to observe, that the entropy formulation of Plotnikov, on the one hand introduce an admissibility condition that is crucial to have uniqueness at least when we consider stable phases and  on the other hand allow  the solution to satisfy the original  forward--backward differential equation in a very weak way (see Definition \ref{defi3}).  This is necessary since  it is not possible to have existence for the classical backward parabolic equation with generic initial data. The main result of this paper is to show that uniqueness of the entropy solution fails in  the general case.  

\noindent Examples of explicit entropy solutions of  the forward--backward parabolic equation that takes values also in the unstable phase are given in \cite{GT}, where it  is considered the ``Riemann problem'' and   a  solution is obtained by self similar methods. More recently in \cite {T2},   it was studied the ``two--phase problem''  where one of the  two  phases is the unstable one.

This paper is organized  in two further sections.

\noindent   In Section \ref {2} we  shall state the precise definition of entropy solution, in particular we shall recall briefly the considerations that lead to this kind of formulation. Moreover,  we shall give a characterization of the entropy solution, showing that this is related to the monotonicity of the coefficients $\lambda_i$, $i=0,1,2$ respect to the variable $t$. More precisely, the coefficients $\lambda_1$, $\lambda_2$ corresponding to the stable phase  tend to increase and the coefficient $\lambda_0$, corresponding to the unstable phase, tends to decrease.  We will be interested to the situation in which we have only unstable phase at initial time, supposing that a solution of the type (\ref{lambda}) appears at positive time (see Definition \ref{meas_regular}). 

In Section \ref{3} we shall show that, choosing properly the initial data, we  obtain infinite solutions that satisfy the entropy formulation of Plotnikov.  This give a negative   answer to  the open question about  the well--posedness of the entropy formulation for general initial data. It is interesting to note, that existence of entropy solutions that have the structure given in \eqref{lambda} are related to the following  inverse parabolic problem  
  \be\label{F1}
 \left\{\begin{array}{lll} &u_t= u_{xx}+f(x) & \hbox{ in } Q_T:=(0,L)\times (0,T)\\
& u_x(0,t)\equiv u_x(L,t)\equiv 0 & \hbox{ in } (0,T)\\
&  u(x,0)=u_0(x)  & \hbox{ in } (0,L)\\
&  u(x,T)=g(x)  & \hbox{ in } (0,L),\end{array}
  \right.
 \ee

\noindent where $u_0$, $g$ are given functions and $f$ is the unknown to be determined. There is a wide literature about this kind of problem, here we just mention  the classic book of Isakov \cite{Is} and  we underline that since we have freedom in the  choice of the  data $u_0$ and $g$ we can easily exhibit  solutions using very classical methods.
 
 \section{Entropy formulation}\label{2}

In this section we  go back quickly to the motivation that are beyond the definition of entropy solution of a  forward--backward parabolic  problem given in \cite {Pl1}, \cite {Pl2}, \cite {Pl3}.
 Then, we characterize the admissibility condition in order to built the counterexample of not uniqueness in Section \ref {3}. 

As said in the Introduction, the idea is to consider the solution as that obtained as singular limit of the approximation problem (\ref{E2}). Problem (\ref{E2}) is analyzed in \cite{NP} also in the multidimensional case. Existence and uniqueness is proved by classical methods of ODE in Banach space (see also \cite{MTT}). Moreover in  \cite{NP}  a viscous entropy inequality is obtained. More precisely, the solution $u_{\epsilon}$ of problem (\ref{E2}) satisfies
\be\label{entri}
\qquad\int\!\!\!\!\int_{Q_T} \!\! \Big \{ G(u^\epsilon)\psi_t  -  g(v^\epsilon) \nabla v^\epsilon \cdot  \nabla \psi
- g'(v^\epsilon)| \nabla v^\epsilon|^2\psi  
\Big \}\,dxdt \, \ge\, 0  \,   
\ee
for any $T>0$, $ \psi \in C^\infty_0(Q_T)$, $ \psi \ge0$, where for any function  $g \in C^1(\R)$, $g' \ge 0$, 
\be\label{defG}
G(u):= \int^u_0 g(\phi(s))ds + K \qquad (K \in \R) \, .
\ee

Using these inequalities and choosing properly the function $g$ it is possible  to obtain a priori estimates in $L^{\infty}$ for $u_{\epsilon}$, $v_{\epsilon}$ that do not depend on $\epsilon$ (see \cite{NP} and \cite{MTT} for the details).

\noindent Using this kind of estimates we can extract proper subsequences $\left\{u_{\epsilon_n}\right\}$, $\left\{v_{\epsilon_n}\right\}$ that converge, respectively,  in the $L^{\infty}$ weak$ ^\star$ topology,  to function $u$ and $v$.  Obviously, we can state that equation $u_t= v_{xx}$ is satisfied in a weak sense but in general $v\neq \phi(u)$ and  we can not  pass to the limit in the viscous entropy inequality (\ref{entri}).

In order to overcome such obstacle, Plotnikov in \cite{Pl2} studied  the Young measure $\nu_{(x,t)}$ associated to the converging sequence $\left\{u_{\epsilon_n}\right\}$, proving that this is {\it is a superposition of Dirac measures} concentrated on the three monotone  branches of the graph of $v=\phi(u)$; the functions $\beta_1$, $\beta_2$, $\beta_0$ defined in  the previous section (see Figure  \ref{funzione2}). 

\noindent More precisely 

\be\label{dirac}
\nu_{(x,t)}(\tau) =\sum_{i=0}^2\lambda_i(x,t)\delta (\tau - \beta_i(v(x,t)))
\ee
 where  $\delta$ is the classical Dirac measure,  $\lambda_i(x,t)\in L^{\infty} (Q_T)$, $\lambda_i\ge 0$, ($i=0,1,2$)   and   $\displaystyle{\sum_{i=0}^2\lambda_i(x,t)=1}$  in $Q_T$.

This implies, that, for every $f\in C(\R)$ 
\be\label{young}
f(u_{\epsilon_n})\stackrel{\!\!*}{\rightharpoonup} \overline f 
 \quad \textrm{ in  } L^{\infty}(Q_T)\, ;
\ee
where
$$
\overline f(x,t):=\int_{\R}f(\tau)\,d\nu_{(x,t)}(\tau)=\sum_{i=0}^2\lambda_i(x,t)f(\beta_i(v(x,t)))
$$
 for a.e. $(x,t) \in Q_T$. Then, choosing $f(u)=\phi(u)$ we deduce the following relation between $u$ and $v$
  
\be\label{eqcost}  u=\lambda_1\beta_1(v)+\lambda_0\beta_0(v)+\lambda_2\beta_2(v).\ee
Moreover, choosing $f(u)=\phi(u)^2$, it easy to prove that $\phi(u_{\epsilon_n})$ converges in the strong topology $L^2$ to the function $v$.

\noindent In fact it is possible to prove  stronger convergences  of the sequence $v_{\epsilon_n}$ to the function $v$ (see \cite {MTT} for details), in particular this is true in the  $L^{2}((0,T), H^1((0,L))$ topology. Using these considerations, we can   pass to the limit along a proper subsequence in the viscous entropy 
inequality (\ref{entri}), proving that (see \cite {Pl2}, \cite {MTT})
the limit couple $(u,v)$ satisfies 
\be \label{entrilim}
\int \!\! \!\int_{Q_T} \!\!\! \Big \{ G^*\psi_t  \!-\!  g(v) \nabla v \cdot  \nabla \psi  - g'(v)| \nabla v|^2\psi 
\Big \}\, dxdt \, \ge \, 0   \,   
\ee
for any $\psi \in C^\infty_0(Q_T)$, $ \psi \ge0$, where 
\be \label{gistar}
G^*(x,t) :=\sum_{i=0}^2\lambda_iG(\beta_i(v(x,t))) \qquad \textrm{ for a.e. } (x,t) \in Q_T \, .
\ee
Then, it is natural to choose (\ref{entrilim}) as the entropy inequality condition for the solution of the forward--backward problem (\ref{B1})
obtained as singular limit of the approximation problem (\ref{E2}). More precisely, we have the following ``natural'' definition.

\begin{definition}\label{defi3} 
An entropy solution  to problem (\ref{B1}) in $Q_T$
is given by  $u, \lambda_0, \lambda_1, \lambda_2 \in L^\infty (Q_T)$, $v\in L^\infty (Q_T) \cap {L^{2}((0,T), H^1(\Omega))}$ such that:
\smallskip

\noindent $(a)$
$\sum\limits_{i=0}^2\lambda_i=1$, $\lambda_i\ge 0$ and there holds:
\be\label{pl1}
u=\sum_{i=0}^2\lambda_i \beta_i(v)\, 
\ee
with $ \lambda_1 =1$ if $v<A$,  $\lambda_2=1$ if $v>B$;
\smallskip

\noindent $(b)$ 
the couple $(u,v)$ is a weak solution of the  equation $u_t=v_{xx}$ in $Q_T$:
\be \label{weak}
\int \!\! \!\int_{Q_T} \!\!\! \Big \{ u\psi_t  \!-\!   v \psi_x  dxdt +\int_{\Omega}u_0(x)\psi(x,0)dx=0  
\ee
for any $\psi\in C^1(\overline{Q}_T),\ \psi(\cdot,T)= 0$ in $\overline{\Omega}$.
\smallskip

\noindent $(c)$ inequality (\ref{entrilim}) is satisfied for any $\psi \in C^\infty_0(Q_T)$, $ \psi \ge0$ and  $g \in C^1(\R)$, $g' \ge 0$. 
\end{definition}

As we anticipate in the Introduction, in general $\phi(u)\neq v$, then condition $b)$ in the previous definition does not imply that the original forward--backward equation is satisfied in a weak way. This will be true if and only if for a.e.  $(x,t)\in Q_T$, one of the  coefficients $\lambda_i(x,t)$ is equal to $1$   and consequently the others are equal to  $0$. Otherwise the equation will be satisfied  in the sense of the measure--valued solution. On the other hand, by construction, there always exists at least one entropy solution in the sense of Definition \ref {defi3} for every initial data in $L^{\infty}.$ This allows to give sense to the backward equation for a general class of initial data.

Obviously there are some natural questions  related to this definition:

\begin{itemize}
\item Can we rewrite entropy condition (\ref{entrilim}) in a more explicit way? In particular which is the consequence of such condition on the coefficients $\lambda_i$?

\item It is possible to state that, in the case in which the initial data takes values only in the two stable phases, the original forward--backward equation $u_t=\phi(u)_{xx}$ is satisfied at least in the distributional sense?

\item Is there uniqueness for the forward--backward problem (\ref{B1}) in the class of entropy solutions introduced  in Definition  \ref {defi3}?

\end{itemize}

Regarding the first question there is the following result obtained in  \cite{Pl2}, \cite{Pl3}.

\begin{theorem}\label{decrease}
Let $(u,v, \lambda_0, \lambda_1, \lambda_2)$ be an entropy solution in the sense of Definition \ref{defi3} to  problem (\ref{B1}) in $Q_T$. Then $\lambda_i(x,\cdot) \in BV_{loc}(0,T)$
for almost every $x\in \Omega$  $(i=0,1,2)$. Moreover,  if  
$$
{\rm ess} \!\!\!\!\!\sup_{t\in (t_1,t_2)} v(x,t)<B \, 
$$ 
for some interval $(t_1,t_2) \subseteq (0,T)$, then $\lambda_1(x,\cdot)$ is not decreasing in $(t_1,t_2)$. 
Similarly, if
$$
{\rm ess} \!\!\!\!\!\inf_{t\in (t_1,t_2)} v(x,t)>A \, 
$$ 
for some interval $(t_1,t_2) \subseteq (0,T)$, then $\lambda_2(x,\cdot)$ is not decreasing in $(t_1,t_2)$. 
\end{theorem}

This result suggests that coefficients $\lambda_1$ and $\lambda_2$ related to the stable phases tend to increase and therefore $\lambda_0$ decrease. In particular   $\lambda_1$ and $\lambda_2$  do not decrease  unless $v=B$ or $v=A$.

In order to give complete answers to the previous questions, at least for a subclass of initial data that takes values in the two stable phases, we introduce  the   ``two--phase problem''.  More precisely, let us consider an initial data $u_0\in L^{\infty}((0,L))$ that satisfies

\begin{equation}\label{H2}
\left\{\begin{array}{ll}
u_0\leq b\  \mbox{in}\ (0, x_0),\\
u_0\geq c\  \mbox{in}\ (x_0, L),\\
\phi(u_0)\in H^1(\Omega).\end{array}\right.
\end{equation}
where $x_0\in (0,L)$.

\noindent In view of the above assumptions $\eqref{H2}$, we look for a solution to problem \eqref{B1} with a particular structure. More precisely, since the initial datum $u_0$ takes values only in the stable phases, we  impose that  solutions to problem \eqref{B1} are again in these phases with  a regular interface separating the rectangle $Q_T$ into two different regions. We require that the unstable phase $(b,c)$ does not influence the dynamics. Then, in accordance with the general entropy formulation given in \cite{Pl2}, the following definition of two--phase solution was done (see \cite{arma}, \cite{ST}).

\smallskip

\begin{definition}\label{def.two.phases}
Let us suppose that $u_0\in L^{\infty}((0,L))$ satisfies \eqref{H2} and $\phi(u_0)\in C((0,L))$. By a two--phase solution to problem \eqref{B1} we mean a triple $(u,v,\xi)$ such that:
\smallskip

\noindent $(i)$ $u\in L^{\infty}(Q_T),\ v\in C(\overline{Q}_T)\cap L^2((0,T);H^1(\Omega))$, and $\xi:[0,T]\to \overline{\Omega}$, $\xi\in C^1([0,T])$, $\xi(0)=x_0$;
\medskip

\noindent $(ii)$ we have:
\begin{equation}\label{eq.fasi}
u=\beta_i(v)\ \ \mbox{in}\ \,V_i \qquad\ (i=1,2)\,,
\end{equation}
where
\begin{eqnarray}\label{eq.V_1}
&&V_1:=\left\{(x,t)\in \overline{Q}_T\,|\ 0< x<\xi(t)\,,\ t\in [0,T]\right\},\\
&&V_2:=\left\{(x,t)\in \overline{Q}_T\,|\ \xi(t)<x<L\,,\ t\in [0,T]\right\},
\end{eqnarray}
and
\begin{equation}\label{eq.interface}
\gamma:=\partial V_1\cap\partial V_2=\left\{(\xi(t),t)\,|\ t\in [0,T]\right\};
\end{equation}

\noindent $(iii)$ $u$ satisfies conditions $b)$ and $c)$ of  Definition \ref{defi3}.

\end{definition}
\smallskip
 
 Obviously equations \eqref{eq.fasi} imply that $v=\phi(u)$ and the coefficients $\lambda_1$,  $\lambda_2$ are, respectively, equal to $I_{V_1}$, $I_{V_2}$, where  $I_E $ denotes  the characteristic function of the set $E$.
 
 \noindent In this class of solutions it is possible to give a characterization of the entropy inequality (\ref{entrilim}) in terms of admissibility condition for the evolution of the  interface $\xi(t)$. More precisely we have the following result (see \cite{ST} for a proof in the more general case)
 
 \begin{proposition} \label {admis}
Let $(u,v,\xi)$ be a two--phase solution of problem \eqref{B1}. Then
 \be\label{admis2}
 \xi'(t)\left\{\begin{array}{ll}
 \le 0&\hbox{ if }v(\xi(t),t)=B\\
 =0&\hbox{ if }v(\xi(t),t)\in(A,B)\\
 \ge 0&\hbox{ if }v(\xi(t),t)=A.\end{array}\right.
 \ee
 \end{proposition}

This means that interface moves only at the critical value $A$, $B$. This is in accordance with the results in Theorem \ref{decrease}; for any fixed value $\overline x$, $\lambda_1(\overline x,\cdot)$ can pass from the value $1$ to the value $0$ ($\lambda_1$ decrease) only at time $\overline t$   when $v(\overline x, \overline t)=B$, analogously $\lambda_2(\overline x,\cdot)$ can pass from the value $1$ to the value $0$ ($\lambda_2$ decrease) only when $v(\overline x, \overline t)=A$. 

\noindent For this class of problems, we can give the answers to all the previous questions. In fact Proposition \ref{admis} gives the characterization requested in the first question. In order to obtain an answer to the second question  in the class  of  data satisfying  condition \eqref{H2}, we have to prove that there is existence of  two--phase entropy solutions introduced in Definition \ref{def.two.phases}. This problem was studied in \cite{arma}, \cite{T} in the piecewise linear case and in \cite{ST} for the nonlinear cubic case. Then local existence was proved in \cite {arma}, \cite{ST}, and global existence in \cite{T}.
Regarding uniqueness in the class of two phase entropy solutions, this is proved in \cite{arma} in the piecewise linear case and in \cite{ST} in the general nonlinear case.   

The purpose of this paper is to prove that uniqueness fails in the general contest of  Definition \ref{defi3}. Then the last question has a negative answer.  The counterexample that we shall give in   Section 3  is  for initial data that take values in the unstable phase. Therefore the question for initial data that takes values only  in the stable phases is still open. However, it is worth to note,  that, the results obtained  for the ``two--phase problem''  suggest a different answer  for this restricted class of initial data. 

In the last part of this section we analyze a new characterization of the entropy condition \eqref{entrilim} for a class of solutions which will be inroduced in Section 3. First of all, observe that, if we impose that one coefficient $\lambda_i$ is equal to $0$, $e.g.$ $\lambda_1\equiv 0$, we can choose one of the other coefficient in function of the third one, $e.g.$ $\lambda_0=1-\lambda_2$. In some sense, we have again a two phase solution but  with a more general structure. 

We have the following

\begin{proposition}\label{new1} Let $v\in C^2(Q_T)$ be such that:


i) there exists $\lambda\in C^1(Q_T)$ such that $\lambda_t\ge 0$, $0\le \lambda\le 1$ in $Q_T$;

ii) $u:=(1-\lambda)\beta_0(v) +\lambda\beta_2(v)$ and 
$u_t=v_{xx}$ in $Q_T.$

Then $v$ satisfies the entropy inequality \eqref{entrilim}.
\end{proposition}

\proof We have to prove that for any $g\in C^1$, $g'\ge0$

\be\label{again}
\int \!\! \!\int_{Q_T} \!\!\! \Big \{ G^*\psi_t  \!-\!  g(v) \nabla v \cdot  \nabla \psi  - g'(v)| \nabla v|^2\psi 
\Big \}\, dxdt \, \ge \, 0   \,   
\ee
for any $\psi \in C^\infty_0(Q_T)$, $ \psi \ge0$, where 
$$
G^*(x,t) :=(1-\lambda) G(\beta_0(v(x,t)))+ \lambda G(\beta_2(v(x,t)))\qquad \textrm{ for a.e. } (x,t) \in Q_T \, .
$$
with $G$ given in \eqref{defG}. 

\noindent Since for the hypothesis $\lambda$ and $v$ are regular functions, we can integrate by parts  in the first member of \eqref{again}, obtaining
$$\int \!\! \!\int_{Q_T} \!\!\! \Big \{- G^*_t\psi  \!+\!  (g(v))_x     \nabla \psi  - g'(v)| v_x|^2\psi 
\Big \}\, dxdt=$$

$$\int \!\! \!\int_{Q_T} \!\!\! \Big (- G^*_t  \!+\!  (g(v))_x \Big ) \psi
\, dxdt.
$$
Then, it is enough to prove that $- G^*_t  \!+\!  (g(v))_x\ge 0$ in $Q_T$. 

\noindent Observe that
$$-G^*_t=\lambda_t G(\beta_0(v(x,t)))-\lambda_t G(\beta_2(v(x,t)))-(1-\lambda)[G(\beta_0(v(x,t)))]_t-\lambda[G(\beta_2(v(x,t)))]_t=$$
$$\lambda_t G(\beta_0(v(x,t)))-\lambda_t G(\beta_2(v(x,t)))-g(v)\left((1-\lambda)(\beta_0)_t+\lambda(\beta_2)_t\right).
$$
Here we use the definition of $G$, that gives 
$$
\left[G(\beta_i(v))\right]_t=g(\phi(\beta_i(v)))[\beta_i(v)]_t =g(v)[\beta_i(v)]_t \,\, \hbox { for i=0,1,2}.
$$
For condition $ii)$ we have
$$
-G^*_t=\lambda_t G(\beta_0(v(x,t)))-\lambda_t G(\beta_2(v(x,t)))-u_t g(v)+[-\lambda_t\beta_0(v)+\lambda_t\beta_2(v)]g(v).$$
Since $u_t=v_{xx}$, we obtain

$$
g(v)v_{xx}-G^*_t=\lambda_t\left[ G(\beta_0(v))-G(\beta_2(v))+(\beta_2(v)-\beta_0(v))g(v)\right].
$$
Then, since $\lambda_t\ge 0$,  it remains to prove that $G(\beta_0(v))-G(\beta_2(v))+(\beta_2(v)-\beta_0(v))g(v)\ge 0$.

Using again the definition of $G$, we get
$$
G(\beta_0(v))-G(\beta_2(v))+(\beta_2(v)-\beta_0(v))g(v)=\int_{\beta_0(v)}^{\beta_2(v)}[g(v)-g(\phi(s))]\,ds
$$
and we obtain the thesis since $g'\ge 0$,  $\beta_0(v)\le \beta_2(v)$ and $v\ge \phi(s)$ for every $s\in (\beta_0(v), \beta_2(v))$.

\qed

Proposition \ref{new1} suggest a way to obtain an entropy solution with unstable phase. Here we have only one coefficient $\lambda$ that correspond to $\lambda_2$. Again the entropy condition is equivalent to the request that the coefficient corresponding to the stable phase does not decrease.

\section{Non--uniqueness results}\label{3}
In order to produce the non existence counterexample,  we consider an entropy solution that  satisfies  the hypothesis of Proposition \ref{new1}. We call this kind of solution a ``two phase, measure--valued, regular solution''  of problem \eqref{B1}. More precisely

\begin{definition}\label{meas_regular} We say that the triple $u\in C^2(Q_T)\cap C((0,L)\times [0,T))$, $\lambda \in C^1(Q_T)\cap C((0,L)\times [0,T))$, $v\in C^2(Q_T)$  is a ``two phase, measure--valued, regular solution''  of problem \eqref{B1} if and only if:

i) $u(x,0)=u_0(x)$ for every $x\in (0,L)$;

ii) there exist $\displaystyle{\lim_{x\to 0^+}v_x(x,t)}=\lim_{x\to L^-}v_x(x,t)=0$ for every $t\in (0,T)$;

iii) $v\ge A$ and $\lambda=1$ if $v>B$;

iv) the functions  $u,\,\lambda,\, v$ are related by following equations:
$$u=(1-\lambda)\beta_0(v) +\lambda\beta_2(v)   \hbox{\,\, in } Q_T,$$
$$u_t=v_{xx}\hbox{\,\, in } Q_T;$$

v) $\lambda_t\ge 0$, $0\le \lambda\le 1$ in $Q_T$.
\end{definition}
Then we have the following
\begin{proposition} A ``two phase, measure--valued, regular solution'' of problem \eqref{B1} is also a solution  in the sense of Definition \ref{defi3} of problem \eqref{B1} 
\end{proposition}
\proof This is consequence of Proposition \ref{new1} that assures that the entropy inequality \eqref{entrilim} is satisfied. The other conditions requested in Definition \ref{defi3} are immediate. 

\qed

\noindent It is natural to consider solutions in which at least at initial time $t=0$ there is not superposition of phases. This means $v(x,0)=\phi (u_0(x))$. Then  $\lambda(x,0)\equiv 1$ in $(0,L)$ or $\lambda(x,0)\equiv 0$ in $(0,L)$. The former case consists of initial data $u_0$ that are  in the stable phase, in this situation $\lambda \equiv 1$ in $Q_T$ and we  obtain a classical solution
solving a forward parabolic problem. In the second case the initial data is in the unstable phase, then we have a purely backward parabolic problem that we can solve by using the auxiliary function $\lambda$. Observe that condition $v)$ of  Definition \ref{meas_regular}   suggests  that at positive  time a superposition of phases appears and in particular the stable phase becomes dominant respect to the unstable one. 

\noindent In the following we always assume that  $\lambda(x,0)\equiv 0$ in $(0,L)$.

Let us fix  $T>0$, $b',\,c'\in (b,c)$ such that $b'<c'$. We choose $g(x)\in C^1 ([0,L])$ such that $g'(0)=g'(L)=0$, $g(x)\in (b',c')$ for any $x\in [0,L]$.

\noindent Let us consider the following problem 
\begin{equation}\label{new2}
 \left\{\begin{array}{lll} &u_t= \phi(u)_{xx} & \hbox{ in } Q_T=\Omega\times (0,T)\\
& u_x(0,t)\equiv u_x(L,t)\equiv 0 & \hbox{ in } (0,T)\\
&  u(x,T)=g(x) & \hbox{ in } \Omega=(0,L).\end{array}
  \right.
\end{equation} 
This is a well--posed parabolic problem, since it is a backward problem with  condition at final time $T$. 

\noindent Denote with  $\overline u$  the unique solution of problem \eqref{new2}. Therefore   for the maximum principle, we have   $\overline u (x,t)\in(b',c')$ for any  $(x,t)\in Q_T$.  Let   $u_0(x):=\overline u(x,0)$. Obviously $\overline u$ is a ``two phase, measure--valued, regular solution'' of problem \eqref {B1}.  We put $\overline v=\phi (\overline u)$, observing that it satisfies the following problem 
\begin{equation}\label{new3}
 \left\{\begin{array}{lll} &(\beta_0(v))_t= v_{xx} & \hbox{ in } Q_T=\Omega\times (0,T)\\
& v_x(0,t)\equiv v_x(L,t)\equiv 0 & \hbox{ in } (0,T)\\
&  v(x,T)=\phi(g(x)) & \hbox{ in } \Omega=(0,L).\end{array}
  \right.
\end{equation} 
Moreover, since $\phi$ is piecewise linear, we have $\overline v(x,0)=\phi(u_0)$.

\noindent Now we want to obtain a different ``two phase, measure--valued, regular solution'' that has the same initial data $u_0$. 

\noindent Let us impose that the triple of function $u,\,\lambda,\, v$ satisfies condition $iv)$ in Definition \ref{meas_regular}.
Then we get
\be\label{schif}[(1-\lambda)\beta_0(v) +\lambda\beta_2(v)  ]_t=v_{xx}.\ee
Let us integrate \eqref{schif} in $(0,t)$, for any fixed $(x,t)\in Q_T$. Since $\lambda(\cdot,0)\equiv 0$, we obtain 
\be\label{schif2}
[(1-\lambda)\beta_0(v) +\lambda\beta_2(v)](x,t)-\beta_0(v(x,0))=\int_0^t v_{xx}(x,s)\,ds. 
\ee

Observe that $\beta_2(v(x,t))=\beta_0(v(x,t))$ if and only if $v(x,t)=A$. In these points  \eqref {schif2} becomes 

\be\label{schif3}\int_0^t v_{xx}(x,s)-[\beta_0(v(x,s))]_s\,ds=0.\ee
On the other hand, in any point $(x,t)\in Q_T$ such that equation \eqref{schif3} is satisfied, we have $\lambda(x,t)=0$ or $v(x,t)=A.$

\noindent In the following we assume $v>A$ in $Q_T$ and we introduce the function $m(x,t)=v_{xx}(x,t)-[\beta_0(v(x,t))]_t$.

\noindent Therefore, we get
\be\label{lamba}\begin{array}{ll}\lambda(x,t)&=\displaystyle{\frac{\beta_0(v(x,0))-\beta_0(v(x,t))+\int_0^t v_{xx}(x,s)\,ds}{\beta_2(v(x,t))-\beta_0(v(x,t))}}=\\
&\displaystyle{\frac{\int_0^t v_{xx}(x,s)-[\beta_0(v(x,s))]_s\,ds}{\beta_2(v(x,t))-\beta_0(v(x,t))}}=\displaystyle{\frac{\int_0^t m(x,s)\,ds}{\beta_2(v(x,t))-\beta_0(v(x,t))}}\end{array}
\ee
 In order to impose the monotonicity condition for the coefficient $\lambda$ we  derivate equation \eqref{lamba} respect to the variable $t$.
 
 Then, we obtain
 \be\label{lamba2}
 \begin{array}{ll}&\lambda_t(x,t)=\displaystyle{\frac{m(x,t)}{\beta_2(v(x,t))-\beta_0(v(x,t))}}-\\&\displaystyle{\frac{[\beta_2(v(x,t))-\beta_0(v(x,t))]_t\int_0^t m(x,s)\,ds}{(\beta_2(v(x,t))-\beta_0(v(x,t)))^2}}.
 \end{array}
 \ee

 It is clear that the sign of $m$ it is strictly related to the monotonicity condition of the coefficient $\lambda$. 

\noindent We have the following result that gives sufficient conditions.

\begin{proposition}\label{prep}
 Suppose that the function $v$  fulfills the following conditions
\begin{equation}\label{new4}
  \!\! \hbox {there exist }  T'\in (0,T), c_1>\!0,\hbox{ s.t. }\,\beta_2(v(x,t))-\beta_0(v(x,t))\ge c_1 \hbox{ in }\, Q_{T'},\ee
 \begin{equation}\label{new5} \!\!  \hbox {there exist }   T''\in (0,T), \,c_2>0, \hbox{ s.t. } m\ge c_2\hbox{ in } Q_{T''},
\ee
 then  there exists $\overline T\in (0,T]$ 
     such that $\lambda_t$ is not negative in $Q_{\overline T}$.
\end{proposition}

The proof  is an immediate consequence of \eqref{lamba2}.
  
   Let us  choose $v(x,t)=\overline v+t$, where $\overline v$ is the solution of problem \eqref{new3}, then function $v$ satisfies the following
\begin{equation}\label{new6}
 \left\{\begin{array}{lll} & v_{xx} -(\beta_0(v))_t=|\sigma|& \hbox{ in } Q_T=\Omega\times (0,T)\\
& v_x(0,t)\equiv v_x(L,t)\equiv 0 & \hbox{ in } (0,T)\\
&  v(x,0)=\phi(u_0) & \hbox{ in } \Omega=(0,L).\end{array}
  \right.
\end{equation} 
here $\sigma=\frac{c-b}{A-B}$ is obtained by the definition of the piecewise function $\phi$ given in \eqref{philin}.  In particular the hypothesis  \eqref{new5} is  satisfied by the function $v$.

Moreover, we can choose $T'$ small enough, such that $v\in (b',c'')$ in $Q_{T'}$ with  $c''\in (c',c)$, then we obtain \eqref {new4} with a proper constant $c_1$. 

\noindent Finally, since $\lambda(\cdot,0)\equiv 0$, using  \eqref{lamba2} we can choose   time $\overline T$   such that  $\lambda \in [0,1)$ for any  $(x,t)\in Q_{\overline T}$. Then we can exhibit  two different ``two phase, measure--valued, regular solution'' of problem \eqref{B1} with  initial data $u_0$. The first one is given by the triple $\overline u$, $\overline \lambda\equiv 0$ and   $\overline v$, and the second one is given by $u,\,\lambda,\, v$ where $v$ is previously defined, $\lambda$ is given  in  \eqref{lamba} and  $u=(1-\lambda)\beta_0(v)+\lambda\beta_2(v)$.

\noindent Let us observe that we can obtain an infinite family of ``two--phase, measure--valued, regular solution'' with $u_0$ as initial condition. Actually we only need to check that conditions \eqref{new4}, \eqref{new5} are satisfied. For example we can find 
solutions of $|\sigma| v_t+v_{xx}=f(x)\ge c_2>0$ that fulfill initial and boundary condition  by standard method of eigenfunction expansion. By straightforward calculation we obtain a solution for every source function   $f\ge c_2>0$, such that
$$
f(x)=\sum_{k=0}^N a_k\cos\left(\frac{k\pi x}{L}\right)
$$
 with  $N\in \N$ and $a_k\in \R$, $k=0\cdots, N$.
 
 Analogous techniques could be used to prove existence of a ``two--phase, measure--valued, regular solution'' with general initial data $u_0$ that takes values in the unstable phase. We are not interested  to consider in detail existence problems. We limit ourself  to highlight  that, in order to obtain existence using Proposition \ref{prep},   it is useful to consider the following parabolic  backward inverse problem, where the unknown data $f$ has to be strictly positive.
 \begin{equation}\label{new7}
 \left\{\begin{array}{lll} & |\sigma|v_t+v_{xx} =f(x,t)& \hbox{ in } Q_T=\Omega\times (0,T)\\
& v_x(0,t)\equiv v_x(L,t)\equiv 0 & \hbox{ in } (0,T)\\
&  v(x,0)=\phi(u_0) & \hbox{ in } \Omega=(0,L).\end{array}
  \right.
\end{equation} 

In general problem \eqref{new7} is underdetermined  unless we fix a final data $v(x,T)=v_T(x)$. On the other hand using the methods of  eigenfunction expansion, we see that it is necessary  to impose some restriction on the initial data. In order to give the idea,  we consider   the simple case in which we suppose that the  solution $f$   depends only on the variable $x$. Then,  if 

$$v_0(x)=\phi(u_0)=\sum_{k=0}^{\infty} a_k\cos\left(\frac{k\pi x}{L}\right)$$
and 
$$v_T(x)=\sum_{k=0}^{\infty} b_k\cos\left(\frac{k\pi x}{L}\right),$$
we obtain 
$$
f(x)=\sum_{k=0}^{\infty} f_k\cos\left(\frac{k\pi x}{L}\right)$$
such that 
$$\!\!\!\!\!\!\!\!\!\!\!\!\!\!\!\!\!\!\!\!\! f_0=\frac{(b_0-a_0)|\sigma|}{T};$$
$$f_k=\frac{\pi^2k^2\left[b_k-a_ke^{\frac{\pi^2Tk^2}{L^2|\sigma|}}\right]}{L^2\left(e^{\frac{\pi^2Tk^2}{L^2|\sigma|}}-1\right)}\quad\,\,k\ge 1.
 $$
This suggest that we can choose properly $v_T$ in order to have a solution $f(x)\ge c_2$  but it is necessary to impose a summability condition to the coefficients $a_k$.


\address{Dipartimento di Matematica ``G. Castelnuovo''\\
Universit\`a di Roma ``La Sapienza''\\
Piazzale A. Moro, 2 --- 00185 Roma, Italy\\

\email{terracin@mat.uniroma1.it}\\

\end{document}

 \ee